\documentclass[11pt]{article}

\usepackage{amsmath,amssymb,epsfig}

\newcommand\zero{{\bf 0}}

\newcommand\bfv{{\bf v}}
\newcommand\x{{\bf x}}
\newcommand\z{{\bf z}}

\newcommand{\cc}{{\mathbb C}}

\newcommand{\bfxi}{{\mbox{\boldmath $\xi$}}}

\newtheorem{theorem}{Theorem}[section]

\newtheorem{proposition}[theorem]{Proposition}
\newtheorem{algo}[theorem]{Algorithm}

\newcommand\qed{{\hspace*{\fill}$\Box$\vskip12pt plus 1pt}}

\begin{document}

\title{Sampling Algebraic Sets in Local Intrinsic Coordinates\thanks{This 
material is based upon work supported by the National Science Foundation
under Grant No.\ 0713018.}}

\author{Yun Guan\thanks{Department of Mathematics, Statistics,
and Computer Science,
University of Illinois at Chicago,
851 South Morgan (M/C 249),
Chicago, IL 60607-7045, USA.
{\tt email:} guan@math.uic.edu
{\tt URL:} http://www.math.uic.edu/{\~{}}guan}
 \and Jan Verschelde\thanks{Department of Mathematics, Statistics,
and Computer Science,
University of Illinois at Chicago,
851 South Morgan (M/C 249),
Chicago, IL 60607-7045, USA.
{\tt email:} jan@math.uic.edu
{\tt URL:} http://www.math.uic.edu/{\~{}}jan}
}

\date{14 December 2009}

\maketitle

\begin{abstract}
Numerical data structures for positive dimensional solution sets of
polynomial systems are sets of generic points cut out by random planes
of complimentary dimension.
We may represent the linear spaces defined by those planes either 
by explicit linear equations or in parametric form.
These descriptions are respectively called extrinsic and intrinsic 
representations.  While intrinsic representations lower the cost of
the linear algebra operations, we observe worse condition numbers.
In this paper we describe the local adaptation of intrinsic coordinates
to improve the numerical conditioning of sampling algebraic sets.
Local intrinsic coordinates also lead to a better stepsize control.
We illustrate our results with Maple experiments
and computations with PHCpack on some benchmark polynomial systems.

\bigskip



\end{abstract}

\section{Motivation, Definitions, and Problem Statement}

A polynomial system $f(\x) = \zero$, $\x = (x_1,x_2,\ldots,x_n)$,
defines an algebraic set $f^{-1}(\zero) \subset \cc^n$.
The polynomials of $f$ belong to~$\cc[\x]$.
We assume (for simplicity of exposition throughout the paper):
\begin{enumerate}
\item $f^{-1}(\zero)$ is pure dimensional, $k$ is its codimension,
      so $\dim(f^{-1}(\zero)) = n-k$;
\item $f(\x) = \zero$ is a complete intersection,
      and in particular: $f = (f_1,f_2,\ldots,f_k)$;
\item $f^{-1}(\zero)$ is reduced, i.e.: of multiplicity one.
\end{enumerate}
To remove the third assumption, a deflation operator~\cite{LVZ06}
(see also~\cite{DZ05})
as proposed in~\cite[\S13.3.2]{SW05} should be applied.
The first two assumptions are made for notational convenience.

The numerical treatment of positive dimensional algebraic sets
was first proposed in~\cite{SW96} and elaborated in a series of
papers by the authors of~\cite{SW05} and the second author,
see also~\cite{SVW9} for another introduction.
The algorithms in numerical algebraic geometry are
implemented in PHCpack~\cite{Ver99} and Bertini~\cite{Bertini}
(see~\cite{BHSW08b} and~\cite{SVW7}) and can be executed via
MATLAB (or Octave)~\cite{GV08b}, Maple~\cite{LV06}, 
and Macaulay~2~\cite{Ley09}.

One of our benchmark examples is a family of systems,
defined by all adjacent minors
of a general 2-by-3 matrix (\cite{DES98}, \cite{HS00}):

\begin{equation}
  \left[
     \begin{array}{cccc}
        x_{11} & x_{12} & x_{13} \\
        x_{21} & x_{22} & x_{23}
     \end{array}
  \right]
  \quad
  f(\x) =
  \left\{
     \begin{array}{rcl}
        x_{11} x_{22} - x_{21} x_{12} \!\!\! & = & \!\!\! 0 \\
        x_{12} x_{23} - x_{22} x_{13} \!\!\! & = & \!\!\! 0.
     \end{array}
  \right.
\end{equation}
For this example, we have $n = 6$, $k = 2$,
and we have a complete intersection: $\dim(f^{-1}(\zero)) = n - k = 4$.
To compute $\deg(f^{-1}(\zero))$, we add $n-k$ general
linear equations $L(\x) = \zero$ to $f(\x) = \zero$
and solve $\{ f(\x) = \zero, L(\x) = \zero \}$.
Generic points on the solution set defined by the system
for all adjacent minors of a general 2-by-3 matrix
satisfy (for random coefficients $c_{ij} \in \cc$):

\begin{equation}
  \left\{
     \begin{array}{rcl}
        x_{11} x_{22} - x_{21} x_{12} \!\!\! & = & \!\!\! 0 \\
        x_{12} x_{23} - x_{22} x_{13} \!\!\! & = & \!\!\! 0 \\
       c_{10} + c_{11} x_{11} + c_{12} x_{12} + c_{13} x_{13}
     + c_{14} x_{21} + c_{15} x_{22} + c_{16} x_{23} \!\!\! & = & \!\!\! 0 \\
       c_{20} + c_{21} x_{11} + c_{22} x_{12} + c_{23} x_{13}
     + c_{24} x_{21} + c_{25} x_{22} + c_{26} x_{23} \!\!\! & = & \!\!\! 0 \\
       c_{30} + c_{31} x_{11} + c_{32} x_{12} + c_{33} x_{13}
     + c_{34} x_{21} + c_{35} x_{22} + c_{36} x_{23} \!\!\! & = & \!\!\! 0 \\
       c_{40} + c_{41} x_{11} + c_{42} x_{12} + c_{43} x_{13}
     + c_{44} x_{21} + c_{45} x_{22} + c_{46} x_{23} \!\!\! & = & \!\!\! 0.
     \end{array}
  \right.
\end{equation}
Except for an algebraic set in the coefficient space $c_{ij}$ for $L$,
the system above has four solutions, we have
four generic points for all adjacent minors
of a general 2-by-3 matrix, so $\deg(f^{-1}(\zero)) = 4$.

To save work, reducing the number of variables from 6 to 2,
we choose a different representation for the
linear space defined by the equations~$L(\x) = \zero$,
representing the 2-plane $L^{-1}(\zero)$ in $\cc^6$ as
\begin{equation}
   \left[
     \begin{array}{c}
        x_{11} \\
        x_{12} \\
        x_{13} \\
        x_{21} \\
        x_{22} \\
        x_{23} 
     \end{array}
   \right]
   = 
   \left[
     \begin{array}{c}
        b_1 \\
        b_2 \\
        b_3 \\
        b_4 \\
        b_5 \\
        b_6
     \end{array}
   \right]
   + \xi_1
   \left[
     \begin{array}{c}
        v_{11} \\
        v_{12} \\
        v_{13} \\
        v_{14} \\
        v_{15} \\
        v_{16}
     \end{array}
   \right]
   + \xi_2
   \left[
     \begin{array}{c}
        v_{21} \\
        v_{22} \\
        v_{23} \\
        v_{24} \\
        v_{25} \\
        v_{26}
     \end{array}
   \right]
\end{equation}
spanned by an offset point ${\bf b} \in \cc^6$
and an orthonormal basis $\{ {\bf v}_1, {\bf v}_2 \}$.
The tuple $(\xi_1, \xi_2)$
defines {\em intrinsic coordinates} for the generic points,
introduced in~\cite{SVW11} to speedup the algorithms of~\cite{SVW10}.

The reduction from six to two variables reduces the cost of solving
linear systems by a factor of nine.  This reduction improves the
efficiency of Newton's method when computing sample points
on the algebraic set, one of the basic operations in
numerical algebraic geometry~\cite{SW05}.

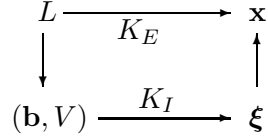
\begin{figure}[hbt]
\begin{center}
\begin{picture}(100,50)
\put(10,40){$L$}
\put(12,35){\vector(0,-1){20}}
\put(0,0){$({\bf b},V)$}
\put(20,43){\vector(1,0){63}}  \put(40,33){$K_E$}
\put(33,3){\vector(1,0){50}}   \put(48,7){$K_I$}
\put(90,40){$\x$}
\put(93,15){\vector(0,+1){20}}
\put(90,0){$\bfxi$}
\end{picture}
\caption{A commutative diagram for extrinsic~$\x$ and intrinsic~$\bfxi$
coordinates of generic points.  The vertical arrows require linear algebra
while the horizontal arrows involve the solution of polynomial systems.
Their sensitivities are determined by condition numbers~$K_E$ and $K_I$.}
\label{figcommutediagram}
\end{center}
\end{figure}

For any $f^{-1}(\zero) \in \cc^n$ with $\dim(f^{-1}(\zero)) = n-k$,
we use a general $k$-plane~$L$ to compute generic points.
This general $k$-plane~$L$ may be defined in two equivalent ways:
\begin{enumerate}
\item $L(\x) = 0$ is a system of $n-k$ general linear equations in $\x$,
\item ${\bf b} \in \cc^n$ is an offset point, and
      $V = [ \bfv_1 ~ \bfv_2 ~ \cdots ~ \bfv_k ] \in \cc^{n \times k}$,
      with $V^* V = I_k$, i.e.: $V$ is an 
      orthonormal\footnote{Although it suffices
to require that the columns of the matrix~$V$ are linearly independent,
the orthonormality condition $V^*V = I_k$
(using complex conjugated inner products and
 $I_k$ is the $k$-by-$k$ identity matrix)
is beneficial.} basis of vectors.
\end{enumerate}
If linear equations define~$L$, solving 
$\{ f(\x) = \zero, L(\x) = \zero \}$ gives generic points in their usual
form that we call an extrinsic coordinate representation.
Using $({\bf b}, V)$ for $L$ gives
intrinsic coordinates $\bfxi = (\xi_1,\xi_2,\ldots,\xi_k)$ for
generic points $\x$:
\begin{equation}
   \x = {\bf b} + \xi_1 \bfv_1 + \xi_2 \bfv_2 + \cdots + \xi_k \bfv_k
      = {\bf b} + V \bfxi.
\end{equation}
With intrinsic coordinates for generic points, the original variables $\x$
become place holders when solving
$f(\x = {\bf b} + V \bfxi) = \zero$.
In Figure~\ref{figcommutediagram}, we outline the two ways to
compute generic points.

In shorthand notation, the general $k$-plane $L$ is represented 
as $({\bf b}, V)$ and we use intrinsic coordinates $\bfxi \in \cc^k$
to denote the generic points.  
When sampling points, the
moving $L$ from $({\bf b}, V)$ to $({\bf c}, W)$,
is done via the obvious homotopy:
\begin{equation} \label{eqobvious}
  f\left(
      \begin{array}{ccccc}
         \x = & (1-t){\bf b} + t {\bf c} & + & 
                ( (1-t) V + t W ) & \bfxi \\
       & \mbox{moving offset point}  &
       & \mbox{moving basis vectors} &
      \end{array}
   \right) = \zero.
\end{equation}
As $t$ moves from 0 to 1, 
the solution paths $\bfxi(t)$ are tracked
with predictor-corrector methods 
and give new generic points on~$f^{-1}(\zero)$.
For introductions to path following and continuation methods
we refer to~\cite{AG03} and~\cite{Mor87}, see also~\cite{Li03}.

While the diagram in Figure~\ref{figcommutediagram} commutes
for exact operations, using floating-point arithmetic forces us
to take into account condition numbers.  
These condition numbers bound the growth of
the relative errors on the solutions as a consequence of relative
errors on the input data.  As we keep~$f$ fixed during the computation
we only consider relative errors on the representations of the $k$-plane~$L$.
Formally, we introduce condition numbers $K_E$ and $K_I$
on the extrinsic and intrinsic coordinate representations
respectively as
\begin{equation}
\frac{||\Delta \x||}{||\x||} \leq K_E \frac{|| \Delta L ||}{|| L ||}
\quad {\rm and} \quad
\frac{||\Delta \bfxi||}{||\bfxi||}
    \leq K_I \frac{|| \Delta ({\bf b},V) ||}{|| ({\bf b}, V) ||}.
\end{equation}
Going to intrinsic coordinates,
we observe a worsening of the numerical conditioning: $K_I \gg K_E$.
Note that the original problem is well conditioned, in other words,
$K_E$ is expected to remain small, because random choices for~$L$
avoids places where the Jacobian matrix of~$f$ drops rank.

To get a first intuition why $K_I \gg K_E$, consider the
binomial expansion of a monomial $x_1^{a_1} x_2^{a_2}$ of~$f$.
If we evaluate $x_1^{a_1} x_2^{a_2}$ 
at $x_1 = b_1 + \xi_1 v_1$ and $x_2 = b_2 + \xi_2 v_2$, we compute
\begin{equation} \label{eqbinexp}
   \left( b_1 + \xi_1 v_1 \right)^{a_1}
   \left( b_2 + \xi_2 v_2 \right)^{a_2} =
     \left( 
       \sum_{i=0}^{a_1} \left( \begin{array}{c} a_1 \\ i \end{array} \right)
       b_1^i (\xi_1 v_1)^{a_1 - i}
     \right)
     \left( 
       \sum_{j=0}^{a_2} \left( \begin{array}{c} a_2 \\ j \end{array} \right)
       b_2^j (\xi_2 v_2)^{a_2 - j}
     \right)
\end{equation}
and we see that any sparse structure of~$f$ will be destroyed.
Moreover, the binomial coefficients in~\eqref{eqbinexp} inflate
the variation among the coefficients in~$f$.

In general, we may write
$f({\bf b} + V (\bfxi + \Delta \bfxi)) =
 f({\bf b} + V \bfxi ) + \Delta f$.
The trouble is that, even for small $|| \Delta \bfxi ||$
we may experience very large $|| \Delta f ||$.

While using multiprecision arithmetic during path tracking~\cite{BHSW08}
may avoid these numerical instabilities, 
using multiprecision numbers significantly slows down the computations
and when the coefficients are known with limited accuracy,
applying multiprecision arithmetic may give misleading answers.
Better stepsize control strategies~\cite{BHSW09} will also be effective
for our problems, but like in dealing with the high powers of the
continuation parameter of polyhedral homotopies~\cite{KK04}, our approach
in this paper is specific to the type of homotopies.
To deal with the numerical instabilities of using intrinsic coordinates,
we propose the use of {\em local} intrinsic coordinates.
We define local coordinates in the next section.
In section~3, we present an algorithm to track a solution path using
intrinsic coordinates, along with an a priori stepsize control evaluation
strategy.  Computational results are discussed in the section~4.

\medskip

\noindent {\bf Acknowledgement.}  We thank Professor Hiroshi Murakami
for his remarks made after the presentation of the first author at the 
session of Symbolic and Numeric Computation at ACA 2009.
His remarks led us to local intrinsic coordinates.

\section{Local Intrinsic Coordinates}

In this section we define local intrinsic coordinate representations
of generic points and address the improved numerical conditioning.

What if we could keep $||\bfxi||$ small?
Writing Greek symbols badly,
the $\xi$ looks close enough to an epsilon,
and then reconsidering the binomial expansions in~(\ref{eqbinexp}):
\begin{eqnarray}
 \left( b_1 + \xi_1 v_1 \right)^{a_1}
    \left( b_2 + \xi_2 v_2 \right)^{a_2} 
 & = & \left( b_1^{a_1} + a_1 b_1^{a_1 - 1} \xi_1 v_1 + O(\xi_1^2) \right)
     \left( b_2^{a_2} + a_2 b_2^{a_2 - 1} \xi_2 v_2 + O(\xi_2^2) \right) \\
 & = & b_1^{a_1} b_2^{a_2} + a_1 b_1^{a_1 - 1} b_2^{a_2} \xi_1 v_1
                         + a_2 b_1^{a_1} b_2^{a_2 - 1} \xi_2 v_2 
     + O(\xi_1^2, \xi_1 \xi_2 , \xi_2^2).
\end{eqnarray}
If we assume that $\xi$ is infinitesimally small,
then we ignore the second order terms $O(\xi_1^2, \xi_1 \xi_2 , \xi_2^2)$.

For general polynomials $f$, 
writing $f({\bf b} + V \bfxi ) = f({\bf b}) + \Delta f$,
the omission of the higher order terms leads to:
$|| \Delta f ||$ is $O(|| V \bfxi ||)$.
Because we may select for the orthonormal basis $V$ 
a nice numerical representation, we have that
$O(|| V \bfxi ||)$ is $O(||\bfxi||)$
and therefore: $|| \Delta f ||$ is $||O(\bfxi)||$.

To keep $||\bfxi||$ small, we now propose to use
the extrinsic coordinates of the generic point
as the offset point for a $k$-plane.
In particular, for $d = \deg(f^{-1}(\zero))$ and $d$ generic points
$\{ \z_1, \z_1 , \ldots, \z_d \}$ on $f^{-1}(\zero)$, consider:
\begin{equation}
   \x = \z_\ell + V \bfxi, \quad \ell = 1,2,\ldots,d.
\end{equation}
Because $L(\z_\ell) = \zero$ for all generic points,
all $\z_\ell + V \bfxi$ represent the same $k$-plane~$L$.
Given an orthonormal basis~$V$ for a $k$-plane
and a set $\{ \z_1, \z_1 , \ldots, \z_d \}$ of $d$
generic points on~$f^{-1}(\zero)$,
the {\em local intrinsic coordinates} to represent $f^{-1}(\zero)$
are defined by the tuple $(\{ \z_1, \z_1 , \ldots, \z_d \}, V)$.

Obviously, the transition from {\em global} intrinsic coordinates
$\x = {\bf b} + V \bfxi$ to {\em local} intrinsic coordinates
is performed by a mere evaluation of~${\bf b} + V \bfxi$.
The close relation between local intrinsic coordinates and
extrinsic coordinates will yield improved condition numbers.

To define the condition number $K_E$ of a zero $\z$ of 
$F := \{ f(\x) = \zero, L(\x) = \zero \}$,
we consider the application of Newton's method:
\begin{equation} \label{eqdefextcond}
    \underbrace{F'(\z)}_{= A} \Delta \z = - F(\z), \quad
    K_E := \kappa(A),
\end{equation}
where $F'$ is the matrix of all partial derivatives of~$F$
and $\kappa(A)$ is the condition number 
of the Jacobian matrix~$A$ of $F$ at $\z$.
Because we assume that~$f(\x) = \zero$ is a complete intersection,
$A \Delta \z = -F(\z)$ is a well defined $n$-by-$n$ linear system.
Strictly speaking, as we keep $f$ fixed and vary only the linear
equations~$L(\x) = \zero$,
we will have $K_E \leq \kappa(A)$, but because generic
points are always well conditioned this distinction is very minor.

In local intrinsic coordinates we replace $\x$ by $\z + V \bfxi$ and
the application of Newton's method leads to
\begin{equation} \label{eqdeflocalcond}
    \underbrace{f'(\z + V\bfxi)}_{= B} \Delta \bfxi 
     = - f(\z + V\bfxi), \quad K_{LI} := \kappa(B),
\end{equation}
where $f'$ is the matrix of all partial derivatives of~$f$
and $\kappa(B)$ is the condition number 
of the Jacobian matrix~$B$ of $f$ at $\bfxi$.
Because we assume that~$f(\x) = \zero$ is a complete intersection,
$B \Delta \bfxi = -f(\bfxi)$ is a well defined $k$-by-$k$ linear system.
We define $\kappa(B)$ as $K_{LI}$,
the condition number of $\z$ represented
in local intrinsic coordinates.

Observe the similarity of~\eqref{eqdefextcond}
with~\eqref{eqdeflocalcond}
as we write~\eqref{eqdefextcond} more explicitly as
\begin{equation}
   \left[
      \begin{array}{c}
         f'(\z) \\ L'(\z)
      \end{array}
   \right]
   \Delta \z =
   - \left[
      \begin{array}{c}
         f(\z) \\ L(\z)
      \end{array}
   \right],
\end{equation}
where $L'$ contains all partial derivatives of~$L$.

\begin{proposition} 
With $K_E$ and $K_{LI}$ 
as defined in~{\rm \eqref{eqdefextcond}}
and~{\rm \eqref{eqdeflocalcond}} respectively: $K_{LI} \approx K_E$.
\end{proposition}
{\em Proof.}  Because in local intrinsic coordinates: $\bfxi = \zero$, 
it does no longer make sense to consider relative errors.  
Moreover, without loss of generality we may always
choose coefficients of the planes so that $||L|| = 1$ and $||V|| = 1$.
Using homogeneous coordinates for~$f$, we assume we work in an appropriate
affine chart so that also~$||\z|| = 1$.  Then the meaning for the
condition numbers $K_E$ and $K_{LI}$ are in the inequalities
\begin{equation} \label{eqcondineqs}
   ||\Delta \z || \leq K_E || \Delta L || \quad {\rm and} \quad
   ||\Delta \bfxi || \leq K_{LI} || \Delta ({\bf b}, V) ||,
\end{equation}
where ${\bf b}$ is an offset point and $V$ are directions
in a parametric representation of~$L$.

Using the commutative diagram of Figure~\ref{figcommutediagram},
we relate $\Delta \z$ and $\Delta \bfxi$:
\begin{equation}
   \z + V(\bfxi + \Delta \bfxi) = \z + \Delta \z
   \quad \Rightarrow \quad
   \Delta \z = V \Delta \bfxi \quad {\rm as}~ \bfxi = \zero.
\end{equation}
Looking at norms: 
\begin{equation}
  || \Delta \z || = || V \Delta \bfxi || = || \Delta \bfxi ||
  \quad {\rm as}~ || V || = 1.
\end{equation}
Because $V^*V = I_k$, all eigenvalues of~$V$ lie on the complex unit
circle and multiplication with~$V$ is norm preserving.

In local intrinsic coordinates, for $\bfxi = \zero$,
changes $\Delta V$ in the orientation of~$L$ do not influence~$\bfxi$.
So we have $\Delta ({\bf b},V) = \Delta {\bf b}$ and
in case we may consider $||\Delta L|| \approx || \Delta ({\bf b},V) ||$. 
Thus in~\eqref{eqcondineqs} we may interchange $K_E$ with $K_{LI}$,
so $K_{LI} \approx K_E$.~\qed

\section{A Rescaling Algorithm}

In this section we consider the sampling of algebraic sets
using local intrinsic coordinates.
We define a rescaling algorithm and address its numerical stability.
In addition, using local intrinsic coordinates leads to a better
stepsize control.

Generic points $\{ \z_1, \z_1 , \ldots, \z_d \}$ are offset points
for a $k$-plane $L$ with directions in the orthonormal matrix $V$.
In local intrinsic coordinates, moving from $(\z_\ell,V)$
to $({\bf b}, W)$, as $t$ goes from~0 to~1, the deformations are defined by

\begin{equation} \label{eqnewhom1}
  f\left( \x = (1-t){\z}_\ell + t {\bf b} + W \bfxi \right) = \zero.
\end{equation}
In contrast to the obvious homotopy in~\eqref{eqobvious},
we see that only the offset point moves.
We immediately switched from the current directions in~$V$
to the new orthonormal basis~$W$ because $\bfxi = \zero$ 
in local intrinsic coordinates.
But this is only a first indication of the potential of working
with local intrinsic coordinates, we can do better than~\eqref{eqnewhom1}.

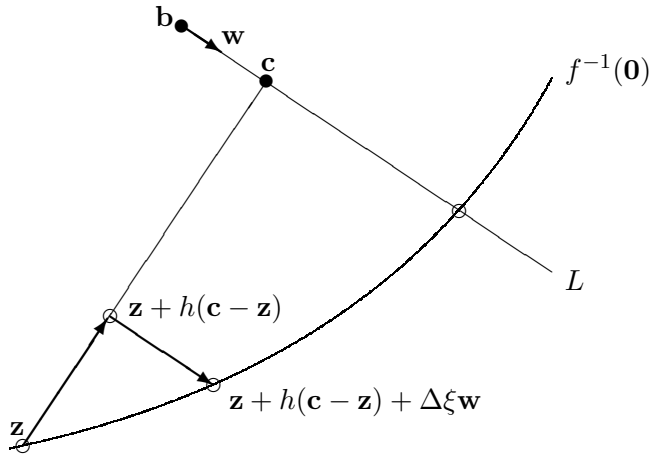
\begin{figure}[hbt]
\begin{center}
\begin{picture}(200,180)(0,0)

\put(0,1){\circle{5}} \put(-5,5){$\z$}
\put(0,1){\line(2,3){92}}
\put(33,50){\circle{5}}
\put(40,50){$\z + h ( {\bf c} - \z )$}

\put(72,24){\circle{5}}
\put(78,15){$\z + h ( {\bf c} - \z ) + \Delta \xi {\bf w}$}

\put(92,139){\circle*{5}} \put(90,143){$\bf c$}

\put(60,160){\circle*{5}} \put(50,160){$\bf b$}
\put(60,160){\line(3,-2){140}} \put(205,60){$L$}

\put(165,90){\circle{5}}

\qbezier(-5,0)(140,30)(200,140) \put(205,140){$f^{-1}(\zero)$}

\thicklines
\put(60,160){\vector(3,-2){15}} \put(75,153){$\bf w$}
\put(33,50){\vector(3,-2){39}}
\put(0,1){\vector(2,3){32}}

\end{picture}
\caption{Schematic of one predictor-corrector step of
the new sampling algorithm, moving from the point~$\z$
to the point where $L$ meets~$f^{-1}(\zero)$.
The line $L$ is defined by ${\bf b} + \xi {\bf w}$.
Using step size $h$, the prediction $h({\bf c} - \z)$ 
added to~$\z$
occurs in a direction orthogonal to~$L$ while the correction
$\Delta \xi {\bf w}$ is parallel to~$L$.  }
\label{figschematic}
\end{center}
\end{figure}

Instead of using~\eqref{eqnewhom1} and moving to $\bf b$,
we point out that any point in the $k$-plane $L$ can serve as an offset point.
Therefore, we should choose the best offset point, 
i.e.: the point closest to the current generic point.  
To compute the closest point,
let $\bf c$ be the orthogonal projection of $\z_\ell$
onto the $k$-plane~$L$.
For some step size~$h$, we then consider:
\begin{equation}
  f\left( \x = {\z}_\ell + h ({\bf c} - \z_\ell) + W \bfxi \right) = \zero
\end{equation}
and apply Newton's method to find the correction $\Delta \bfxi$,
as illustrated in Figure~\ref{figschematic}.

After each step, we add the correction term 
($\Delta \xi {\bf w}$ in Figure~\ref{figschematic})
to the offset point, rescaling the intrinsic coordinates 
to local intrinsic coordinates at the end of the correction stage.
Pseudocode for one predictor-corrector step
is given in Algorithm~\ref{algonestep}, going from one generic point
$\z \in f^{-1}(\zero) \cap K$, where $K$ is the current $k$-plane,
towards~$L$ the target $k$-plane.

\begin{samepage}
\begin{algo} [one predictor-corrector step in local intrinsic coordinates]
\label{algonestep} ~~

{\rm
\noindent~\begin{tabular}{llcr}
\multicolumn{2}{l}{Input: $f = (f_1,f_2,\ldots,f_k), f_i(\x) \in \cc[\x]$,
   $i=1,2,\ldots,k$;}
   & & {\em $\dim(f^{-1}(\zero)) = n-k$} \\
\multicolumn{2}{l}{\hspace{1cm} ${\bf b} \in \cc^n$;}
   & & {\em offset point of $k$-plane $L$} \\
\multicolumn{2}{l}{\hspace{1cm}
     $W = [{\bf w}_1~{\bf w}_2~\cdots~{\bf w}_k] \in \cc^{n \times k}$, 
     $W^*W = I_k$;}
   & & {\em orthonormal basis for $L$} \\
\multicolumn{2}{l}{\hspace{1cm} $\z \in \cc^n$: $f(\z) = \zero$,
     $K(\z) = \zero$;}
   & & {\em generic point on $k$-plane~$K$} \\
\multicolumn{2}{l}{\hspace{1cm} $h > 0$;}
   & & {\em step size} \\
\multicolumn{2}{l}{\hspace{1cm} $\epsilon > 0$.}
   & & {\em accuracy requirement} \\
\multicolumn{2}{l}{Output:
     $\widehat{\z}$, $f(\widehat{\z}) = \zero$, $L(\widehat{\z}) = \zero$:
     $|| \widehat{\z} - {\bf b} || < || \z - {\bf b} ||$.}
   & & {\em generic point closer to~$L$} \\
%
\\
1. & ${\bf v} := \z - {\bf b}$;
   & & {\em go towards offset point} \\
2. & ${\displaystyle {\bf v} := {\bf v}
   - \sum_{i=1}^k ( \overline{{\bf w}_i}^T {\bf v} ) {\bf w}_i}$;
   & & {\em move perpendicular to~$L$} \\
3. & ${\bf v} := \frac{\bf v}{||{\bf v}||}$;
   & & {\em normalize so $|| {\bf v} || = 1$} \\ 
4. & $\widetilde{\z} := \z + h ~ {\bf v}$; 
   & & {\em prediction for new generic point} \\
5. & $\widehat{\z} := \widetilde{\z}$; 
$\bfxi := \zero$; & & {\em initialize for Newton corrector} \\
%
6. & while $||f(\widehat{\z}+W \bfxi)|| > \epsilon$ do 
   & & {\em as long as not accurate enough} \\
6.1 & \hspace{5mm} $\Delta \bfxi 
   := f(\widehat{\z} + W \bfxi)/f'(\widehat{\z} + W \bfxi)$; 
   & & {\em solve a linear system for~$\Delta \bfxi$} \\
6.2 & \hspace{5mm} $\bfxi := \bfxi + \Delta \bfxi$;
   & & {\em update correction} \\
7. & $\widehat{\z} := \widehat{\z} + W \bfxi$.
   & & {\em rescale to local coordinates}
\end{tabular}
}
\end{algo}
\end{samepage}
The orthonormality condition $W^* W = I_k$ is important
for instruction~2 in the Algorithm~\ref{algonestep}
because we can compute the projection just via inner products.
The number of arithmetical operations needed to carry out
instruction~2 in Algorithm~\ref{algonestep} is~$O(kn)$.
Without the condition~$W^*W = I_k$, this cost 
(e.g. via Gram-Schmidt orthogonalization) would 
be at least~$O(k n^2)$.

For the numerical stability of Algorithm~\ref{algonestep},
we first discuss the relationship between the step size~$h$
and the accuracy requirement~$\epsilon$.
If on the one hand $h$ is too small, then the condition
$||f(\widehat{\z}+W \bfxi)|| > \epsilon$ in the while-do
instruction~6 of Algorithm~\ref{algonestep} is directly satisfied.
On the other hand, if $h$ is too large, satisfying the accuracy
requirement of instruction~6 may require too many iterations,
or Newton's method may not converge at all.
We point out that the cost of instruction~6.1 is~$O(k^3)$
and if $f$ is sufficiently sparse (if evaluation and differentiation
go fast), then the cost of execution of Newton's method dominates
the cost of Algorithm~\ref{algonestep}.

In general path tracking algorithms, the step size~$h$ is determined
via a feedback mechanism.  If Newton's method does not converge fast
enough, then the step size is reduced.  If Newton's method needs only
two steps or less, then the step size might be enlarged.
See~\cite{BHSW09} for stepsize control strategies.
The problem with this feedback mechanism is that it comes at the great
expense of the most costly portion of the predictor-corrector method,
i.e.: each reduction of~$h$ comes at the expense of a failed and thus
wasted Newton step.
With local intrinsic coordinates, we can predict the fitness of the
step size with a simple evaluation.
For some step size~$h$ and direction~$\bf v$,
we evaluate and estimate the residual as
\begin{equation}
   ||f( \x = \z_\ell + h {\bf v} )||
    ~{\rm is}~ ||f(\z_\ell) + O(h)|| ~{\rm is}~ O(h).
\end{equation}
For example, if $h = 10^{-2}$ and we see that the residual is
$O(10^{-2})$, then it is fair to expect that after one iteration
of Newton's method, the residual becomes $O(10^{-4})$,
and then $O(10^{-8})$ after the second iteration.

In Algorithm~\ref{algstepcontrol} we define how to cut back on
the step size just by evaluation, {\em before} the start of the
Newton correction.

\begin{algo} [a priori stepsize control by evaluation]
\label{algstepcontrol} ~~

{\rm
\begin{tabular}{llcr}
\multicolumn{2}{l}{Input: $f = (f_1,f_2,\ldots,f_k), f_i(\x) \in \cc[\x]$
   $i=1,2,\ldots,k$;}
   & & {\em $\dim(f^{-1}(\zero)) = n-k$} \\
\multicolumn{2}{l}{\hspace{1cm} $\z \in \cc^n$: $f(\z) = \zero$,
     $K(\z) = \zero$;}
   & & {\em generic point on $k$-plane~$K$} \\
%
\multicolumn{2}{l}{\hspace{1cm} ${\bf v} \in \cc^n$, $||{\bf v}|| = 1$;}
   & & {\em direction vector} \\
\multicolumn{2}{l}{\hspace{1cm} $h > 0$;}
   & & {\em current step size} \\
\multicolumn{2}{l}{\hspace{1cm} $\delta > 0$.}
   & & {\em threshold to reduce $h$} \\
\multicolumn{2}{l}{\hspace{1cm} $1 > \rho > 0$.}
   & & {\em reduction factor for $h$} \\
\multicolumn{2}{l}{Output: $h > 0$.}
   & & {\em updated step size} \\
\\
1. & $y := ||f(\z + h {\bf v})||$; 
   & & {\em evaluate the predicted point} \\
2. & if $y/h > \delta$ then $h := \rho h$.
   & & {\em reduce the step size} 
\end{tabular}
}

\end{algo}
The reduction of the step size in instruction~2 
of Algorithm~\ref{algstepcontrol} could be followed by another 
evaluation of~$f$ to see if $y$ is reduced sufficiently or has
become even too small.

Algorithm~\ref{algstepcontrol} is called after instruction~3
of Algorithm~\ref{algonestep}.

By application of Algorithm~\ref{algstepcontrol},
occurrences of a diverging Newton's method can be greatly reduced
because the size of the residual $||f(\x = \z + W\bfxi)||$ is $O(h)$.

We conclude with a quick cost estimate for the total number of
Newton steps along one path.
In sampling for generic points, we typically choose the new random coefficients
for the $k$-plane as complex numbers on the unit circle, so the distance
between two $k$-planes (and in particular their offset points) is~$O(1)$.
For $h$: $0 < h < 1$, we can see that the total number of Newton iterations
along a solution path is proportional to $1/h$.  For example if $h=0.01$
and we need about~2 or~3 Newton iterations per step, then the total number
of Newton iterations along a solution path will vary between 200 and~300.

The homotopy continuation methods of this paper are different from the
so-called linear homotopies for which an experimental study to certify
path tracking recently appeared in~\cite{BL09}.
A potential future research direction could be to expand the quick
cost estimate of the previous paragraph into a formal complexity study,
along the lines of~\cite{BCSS98} and~\cite{Pet07}.

\section{Computational Results}

Local intrinsic coordinates are available in version 2.3.53
of PHCpack~\cite{Ver99}.  We first describe numerical experiments
done with Maple to compare condition numbers of generic points
on a hypersurface of polynomials of increasing degrees.
Then we report preliminary results on small benchmark problems
with the sampling routines in PHCpack.
All computations were done on one core of a Mac OS X 3.2 Ghz Intel Xeon.

\subsection{Condition Number Estimates}

In this section we illustrate the worsening of the conditioning
of using global intrinsic coordinates on one sparse polynomial.
We give data on sampling with zero and nonzero offset vectors
and relate this experiment to using local intrinsic coordinates.

To estimate the condition numbers we use
{\tt LinearAlgebra[EigenConditionNumbers]} of Maple~12,
with {\tt UseHardwareFloats} set to {\tt true}, see~\cite[Chapter~4]{Mon08}.
The corresponding documentation pages of Maple~12 refer 
to~\cite{And99}.  For an introduction to the perturbation
theory of eigenvalues, see e.g.: \cite[\S4.3]{Dem97}.

We consider one sparse polynomial~$f$ in $n= 10$ variables,
of increasing degrees $d$, starting with $t$ terms.
In addition, we add all the linear terms $c_i x_i$, $i=1,2,\ldots,n$,
to avoid ending up with the origin as a multiple root.
The coefficients are taken on the complex unit circle.
The particular Maple commands used to generate an~$f$ are
\begin{verbatim}
[> n := 10: d := 10: t := 5:
[> c := () -> exp(I*stats[random,uniform[0,2*Pi]](1)):
[> X := [seq(x[i],i=1..n)]:
[> f := X[1]^d + randpoly(X,coeffs=c,degree=d-1,terms=5) + sum(c()*x[i],i=1..n);
\end{verbatim}
The first term of {\tt f} ensures that we have a monic polynomial
after substitution $f({\bf v}\xi)$, for $v_1 = 1$.
That $f$ is monic is convenient for the connection with
the companion matrix when we look at the condition numbers
of the corresponding eigenvalue problem.

To introduce the idea of using different coordinate systems,
we respectively use
\begin{equation}
   \x = {\bf b} + {\bf v} \xi \quad {\rm and} \quad
   \x = {\bf v} \xi, \quad
   {\bf b}, {\bf v} \in \cc^n, 
\end{equation}
where all coefficients in the vectors are also taken on the 
complex unit circle.  With $f({\bf v}\xi) = 0$ we obtain still
a sparse polynomial with all coefficients on the complex unit circle,
which is not the case with $f({\bf b} + {\bf v}\xi) = 0$.
The offset vector of ${\bf b} + {\bf v} \xi$ is responsible for
the variation in the coefficients and the fluctuation of the condition
numbers we observe in our numerical experiments,
summarized in Table~\ref{tabmaplecalc}.

\begin{table}[hbt]
\begin{center}
\begin{tabular}{|c|cc|cc|cc|} \hline
  degrees & \multicolumn{2}{c|}{$f({\bf b} + {\bf v} \xi) = 0$}
          & \multicolumn{2}{c|}{$f({\bf v} \xi) = 0$}
          & ratios of & ratios of \\
  of $f$ & largest & smallest & largest & smallest 
         & smallest & largest \\ \hline
10 & {\tt 5.91e-01} & {\tt 9.02e-02} & {\tt 8.81e-01} & {\tt 4.01e-01}
   & {\tt 6.55e+00} & {\tt 2.20e+00} \\
20 & {\tt 2.77e-01} & {\tt 1.76e-03} & {\tt 8.91e-01} & {\tt 3.31e-01}
   & {\tt 1.57e+02} & {\tt 2.70e+00} \\
30 & {\tt 2.75e-01} & {\tt 6.16e-05} & {\tt 9.49e-01} & {\tt 7.25e-02}
   & {\tt 4.47e+03} & {\tt 1.31e+01} \\
40 & {\tt 4.53e-01} & {\tt 7.14e-06} & {\tt 9.69e-01} & {\tt 1.87e-01}
   & {\tt 6.34e+04} & {\tt 5.17e+00} \\ \hline 
\end{tabular}
\caption{Estimates for
the inverse condition numbers of eigenvalues of the companion matrices
of $f({\bf b} + {\bf v} \xi) = 0$ and $f({\bf v} \xi) = 0$.
For degrees $d = 10$, 20, 30, and 40, we list the largest and smallest
inverse condition numbers.  For $f({\bf v}\xi) = 0$, we see the range
between smallest and largest not widen that much, whereas for
$f({\bf b} + {\bf v}\xi) = 0$, the conditioning steadily worsens.  }
\label{tabmaplecalc}
\end{center}
\end{table}

As we see from Table~\ref{tabmaplecalc}, all roots of $f({\bf v} \xi) = 0$
are well conditioned.  To compare the conditioning of local intrinsic
coordinates, we take the first root $z_1$ of $f({\bf v} \xi) = 0$
and consider the companion matrix $A$ of $f(z_1 + {\bf v} \xi) = 0$.
For increasing degrees,
the condition number for the zero $\xi = 0$ corresponding to $z_1$
is always reported as {\tt 1.00e+00}.
The smallest inverse condition numbers of the eigenvalues of~$A$
for degrees $d = 10$, 20, and 30 are respectively 
{\tt 8.42e-05}, {\tt 1.08e-12}, and {\tt 3.69e-14}.
This implies that for $d = 30$ we have lost all accuracy
as our working precision are the standard hardware floats.

In this simple Maple experiment we illustrate that, although sampling
a hypersurface is reduced to solving univariate polynomial equations,
for hypersurfaces defined by polynomials of high degrees we cannot
use the same representation of a general line to define generic points.
If we adapt the offset point and switch to local intrinsic coordinates,
then the generic points are well conditioned.

\subsection{Sampling Benchmark Systems}

The input to the sampling problem is one set of generic points
on $f^{-1}(\zero) \cap L$ and a new $k$-plane~$K$.
On output is a new set of generic points on~$f^{-1}(\zero) \cap K$.

The polynomial systems we selected occur in the literature.
We briefly summarize the main characteristics of these systems:

\begin{enumerate}
\item All adjacent minors of a general 2-by-$n$ matrix, $n = 3,4,\ldots$.
      This is a family of nice quadratic equations
      arising in algebraic statistics~\cite{DES98}.
\item The cyclic $n$-roots systems are well known academic benchmarks.
      If $n$ has a quadratic divisor, then the system has a positive
      dimensional solution set~\cite{Bac89}.
      In our experiments we use the cyclic 8-roots system,
      which has a one dimensional solution set of degree~144.
\item Griffis-Duffy platforms~\cite{GD93}
      are architecturally singular mechanisms~\cite{HK00}, their
      motion correspond to curves of degree 40 in 8-space~\cite{SVW6}.
\end{enumerate}
For the purposes of this paper, the computation of the first set of
generic points is considered as given, typically in extrinsic
coordinate representation.  

Except for the adjacent minors, the systems are not complete intersections.
For $m > k$, to make an $m$-by-$k$ system $f$ square, 
we generate a random $k$-by-$m$ matrix $C$ and work with~$C \times f$.

To test the improvement from using local intrinsic coordinates,
we sample new generic points from the solution sets.
Our computational experimental setup consists of three stages:
\newline (1) Given one set of generic points, 
             we generate another random $k$-plane~$L$.
\newline (2) We then move the given set of generic points to lie on~$L$.
\newline (3) At the end we check results for accuracy, 
             count \#predictor-corrector steps, record elapsed cpu times.

Note that the recorded cpu times are only meant to give an indication
on the relative practical difficulties of these problems.
More relevant are the number of iterations performed by Newton's method
along the paths.

In Table~\ref{tabexperiments} we summarize the results.
Even as the systems we selected as benchmark examples are not challenging, 
we observe a clear benefit of using local intrinsic coordinates,
even for the systems defined by quadratic equations.
The benefit is perhaps most significant for the cyclic 8-roots problem
where the degree of the $i$th polynomial equals~$i$.

\begin{table}[hbt]
\begin{center}
\begin{tabular}{|c||r|r|r||r|r|} \hline
polynomial system & $n$ & $n-k$ & $d$~~ & \#iterations & timings
 \\ \hline \hline
Griffis-Duffy platform & 8 &  1~~ & 40  & 207/164 & 550/535~$\mu$sec \\ \hline
 cyclic 8-roots system & 8 &  1~~ & 144 & 319/174 & 5.3/3.2 sec \\ \hline
  all adjacent minors & 22 & 12~~ & 1,024 & 285/219 & 44.6/40.3 sec \\ \hline
\end{tabular}
\caption{Preliminary experiments on three systems.
For each system we respectively list the ambient dimension~$n$,
the dimension $n-k$ of the solution set, and the degree $d$ of the set.
We list the average number of Newton iterations along a path
for intrinsic and local intrinsic coordinates, as well 
as user cpu timings. }
\label{tabexperiments}
\end{center}
\end{table}

\section{Conclusions}

We list at least three
advantages of using local intrinsic coordinates for sampling:
(1) only the offset point moves;
(2) the sparse structure of the polynomials is kept; and
(3) we can control the step size by evaluation.
Applications to numerical algebraic geometry include
(1) implicitization via interpolation;
(2) monodromy breakup algorithm; and
(3) diagonal homotopies.
In particular, local intrinsic coordinates will add to the
robustness of our parallel subsystem-by-subsystem solver~\cite{GV08}.

\bibliographystyle{plain}

\end{document}